\documentclass[a4paper,12pt,twoside]{article}
\usepackage{amsfonts,amsmath,amssymb}
\usepackage[body={4.7in,7.5in},top=1.5in,left=1.5in]{geometry}
\numberwithin{figure}{section} \numberwithin{table}{section}

\usepackage{enumerate}
\usepackage{graphicx}
\usepackage{array}
\usepackage{fancyhdr}
\def\refhg{\hangindent=20pt\hangafter=1}
\def\refmark{\par\vskip 1mm\noindent\refhg}

\title{A Generalization of the Marshal-Olkin Scheme and Related Processes}
\author{Satheesh S$^1$, Sandhya E $^2$ and Prasanth C B $^3$\\
\small $^1$Department of Applied Sciences, Vidya Academy of \\\small Science and Technology, Thalakkottukara, Thrissur-680 501, India.\\
\small $^2$Department of Statistics, Prajyoti Niketan College\\\small Pudukkad, Thrissur-680 301, India.\\
\small $^3$Department of Statistics, St. Thomas College, Pala-686 574, India.\\
\small e-mail ids:$^1$ssatheesh1963@yahoo.co.in; $^2$esandhya@hotmail.com \\\small and $^3$cbpwarrier@yahoo.co.in\\
}
\date{}
\numberwithin{equation}{section}

\begin{document}

\maketitle
\chead{}\rhead{}\lfoot{}\cfoot{\thepage}\rfoot{}
\renewcommand{\headrule}{\vbox to 0pt {\hbox to \headwidth{}\vss}}
 
\noindent{\textbf{Abstract}\\ }
      \small
A generalization of the Marshal-Olkin parametrization scheme is developed and stochastic models related to it are discussed here.
\noindent{\textbf{AMS(2010) Subject Classification Numbers:}
60E05, 60E07, 60G70, 62M10.}
\section{Introduction.}

Given a distribution, the Marshal-Olkin parametrization scheme (M-O scheme) gives a generalization of it in terms of its survival function (\textit{s.f}) by adding a parameter to it, thus providing more flexibility, Marshal and Olkin (1997). For a given \textit{s.f} $\overline{F}$, the M-O scheme is described by the \textit{s.f} 
  
\begin{equation} \overline{G}(x,\alpha) = \frac{\alpha \overline{F}(x)}{1-(1-\alpha) \overline{F}(x)}, x\in R, \alpha>0. \end{equation} 

Satheesh and Nair (2004) observed that the M-O scheme has essentially a geometric-minimum structure. Their investigation also suggested a generalization of the M-O scheme as follows. For a given \textit{s.f} $\overline{F}$, this is described, for $k>0$ integer and $a>0$, by the \textit{s.f} 

\begin{equation} \overline{H}(x,a) =  \left\{\frac{\overline{F}^k(x)}{a-(a-1)\overline{F}^k(x)}\right\}^{1/k}, x\in R. \end{equation}

This scheme thus adds two parameters to make it more flexible. When $k=1$, (1.2) reduces to the M-O scheme. This generalization has a Harris-minimum structure since replacing $\overline{F}(x)$  by $s\in(0,1)$ and restricting  $a>1$, the RHS of (1.2) is the probability generating function of a Harris($a,k$) distribution studied in detail by Sandhya \textit{et al}. (2008).\\

With reference to the M-O scheme it is to be noted that if $F$ is used instead of $\overline{F}$ in (1.1) we will get $G(x,1/\alpha)=1-\overline{G}(x,\alpha)$. In other words, the distribution function (\textit{d.f}) corresponding to the \textit{s.f} defined in (1.1) is similar to the \textit{s.f} but for a change in the parameter $\alpha$ as given above. But, this is not true for the generalised schemes (1.2) and (2.1) discussed in this paper. Hence there is relevance in studying these schemes seperately.\\

The purpose of this note is to generalize the M-O scheme on the lines of (1.2) in terms of \textit{d.f} and discuss stochastic models related to it. Possible applications of this model are in reliability studies of parallel systems.  
 
\section{A construction and invariance of the generalized M-O scheme.}
\vskip.3cm 

We first show that (1.2) is indeed a \textit{s.f}. Notice that for any \textit{s.f} $\overline{F}(x)$, $\overline{F}^\nu(x), \nu>0$ is always a \textit{s.f}. Thus $\overline{H}(x,\alpha)= \alpha \overline{F}^k(x)/\{1-(1-\alpha)\overline{F}^k(x)\}$ is a \textit{s.f} by the M-O scheme, with $\nu=k$  integer. Since $\overline{H}^{1/k}(x,\alpha)$ is also a \textit{s.f}, the generalization at (1.2) is a \textit{s.f} where $a=\frac{1}{\alpha}$. Further, this scheme is closed under Harris($b,k$)-minimum ($b>1$), but not under Harris($b,k$)-maximum, Satheesh and Nair (2004).\\

Since for any \textit{d.f}  $F(x)$, $F^\nu(x), \nu>0$ is also a  \textit{d.f} we have the following parametrization scheme for \textit{d.f}s, for $k>0$ integer and $a>0$.

\begin{equation} 
H(x,\alpha) = \left\{\frac{F^k(x)}{a-(a-1)F^k(x)}\right\}^{1/k}, x\in R. \end{equation} 
    
Notice that (2.1) is closed under Harris($b,k$)-maximum ($b>1$), but not under Harris($b,k$)-minimum, Satheesh and Nair (2004).\\
\\
\textbf{Theorem 2.1} If $H(x)=\frac{1}{1+\psi(x)}$ is a \textit{d.f}, then $\left\{\frac{1}{1+a\psi(x)}\right\}^{1/k}, k>0$ integer, $a>0$ is also a \textit{d.f}.\\
\\ 
\textit{Proof}. Since $  \left\{\frac{1}{1+a\psi(x)}\right\}^{1/k}=\frac{\{1+\psi(x)\}^{-1/k}}{\left\{\frac{1+a\psi(x)}{1+\psi(x)}\right\}^{1/k}} 
=\frac{\{1+\psi(x)\}^{-1/k}}{\left\{\frac{1+a\psi(x)+a-a}{1+\psi(x)}\right\}^{1/k}},$ we have \begin{equation} 
\left\{\frac{1}{1+a\psi(x)}\right\}^{1/k}=\left\{\frac{H(x)}{a-(a-1)H(x)}\right\}^{1/k}. 
\end{equation}  
Since for any \textit{d.f} $H(x)$, $H^{1/k}(x)$ is a \textit{d.f}, the RHS in (2.2) is the parametrization (2.1) for $H^{1/k}(x)$, proving the assertion.\\
\\
\textbf{Remark 2.1} This theorem suggests a convenient way to construct \textit{d.f}s that are invariant under Harris($b,k$)-maximum since for any \textit{d.f} $H(x)$ we can take $\psi(x)=-logF(x)$ in $\{x: F(x)>0\}$. See also remark 3.1.\\
\\
\textbf{Corollary 2.1} Since for any \textit{d.f} $H(x)$, $H^k(x)$ is a \textit{d.f}, $\frac{1}{1+a\psi(x)}$ is also a \textit{d.f}.\\
\\ 
\textbf{Corollary 2.2} When $a>1$, we get the Harris($a,k$)-maximum of $H^{1/k}(x)$.\\
\\ 
\textbf{Remark 2.2} In particular, if  $\psi(x)=a\psi(cx)$, for some $0<c<1<a$, then theorem 2.1 shows that the \textit{d.f} $H^{1/k}(x)$ is invariant under Harris($a,k$)-maximum upto a scale change $c$. We may call such \textit{d.f}s Harris-max-semi-stable($a,c,k$).\\

In terms of CFs, theorem 2.1 generalizes lemma 3.1 in Pillai (1990). In general, for a characteristic function (CF) $f$, $f^\nu, \nu >0$ is not a CF, though for an infinitely divisible (ID) CF this is true. Thus for those CFs for which this is true the operation analogous to (2.1) (or (1.2)) on a CF results in a CF that is closed under Harris($b,k$)-sum. If we are considering semi-stable CFs, then they are ID and an analogue of theorem 2.1 implies that the CF $f(t)=\left\{\frac{1}{1+\omega(t)}\right\}^{1/k},  \omega(t)=a\omega(ct)$, for some $0<c<1<a$, is invariant under Harris($a,k$)-sum.

\section{Stochastic processes related to Harris-maximum.}
\vskip.3cm  
Here we discuss two stochastic processes where the parameterization (2.1) or \textit{d.f}s of the form $\left\{\frac{1}{1+a\psi(x)}\right\}^{1/k}, k>0$ integer, $a>0$, appear naturally.\\

Pancheva \textit{et al.} (2006) discussed random time changed or compound extremal process (EP) and their theorem 3.2 together with Property 3.2 reads: Let $\{Y(t), t\geq 0\}$ be an EP with homogeneous max-increments and \textit{d.f} $F_t(y)=exp\{-t\mu[(\lambda,y)^c]\}, y\geq\lambda,\lambda>0$ being the
bottom of the rectangle $\{F > 0\}$ and $\mu$ the exponential measure of $Y(1)$, that is, $\mu[(\lambda,y)^c]= -log F(y)$. Let $\{T(t),t\geq0\}$ be a non-negative process independent of $\{Y(t)\}$
having stationary, independent and additive increments with Laplace transform $\varphi^t$. If $\{X(t), t\geq0\}$ is the compound EP obtained by randomizing the time parameter of $\{Y(t)\}$ by
$\{T(t)\}$ then $X(t)=Y(T(t))$. Its \textit{d.f} is $P\{X(t)\leq x\} = \{\varphi[\mu(\lambda, x)^c]\}^t$ which is $\varphi$-max-ID, Satheesh \textit{et al}. (2008). Pancheva \textit{et al}. (2006) also showed that in this set up $\{Y(T(t))\}$ is also an EP. We now have\\
\\
\textbf{Theorem 3.1} The EP obtained by compounding a homogeneous EP is gamma-max-ID if the compounding process is gamma($\alpha,\beta$).\\
\\ 
\textit{Proof}. If $\{Y(t)\}$ is an EP with homogeneous max-increments and \textit{d.f} $e^{-\xi(x)}$, $\{T(t)\}$ a gamma($\alpha,\beta$) process with stationary, independent and additive increments and \textit{d.f} $G$, then the \textit{d.f} of the process $\{Y(T(t))\}$ is given by  $\int^{\infty}_{0}e^{-t\xi(x)}dG(t)=\left\{\frac{1}{1+\alpha \xi(x)}\right\}^\beta$, proving the assertion.\\
\\ 
\textbf{Remark 3.1} Clearly this is a gamma mixture. When $\beta=\frac1{k}$, this \textit{d.f} is the same as the one in
theorem 2.1 and is closed under Harris($b,k$)-maximum ($b>1$). Further, when $a>1$ this \textit{d.f} is itself a Harris($a,k$)-maximum. \\
 
Another stochastic model is the max-AR(1) process described by \textit{i.i.d r.v}s $\{Y_{i,n}\}$ and innovations (\textit{i.i.d r.v}s) $\{\epsilon_{i,n}, i = 1,2,...,k\}$ for a fixed positive integer $k$. Such generalized models were considered by Satheesh \textit{et al.} (2008).

\begin{equation}
\bigvee ^{k}_{i=1} Y_{i,n} = \begin{cases}
\bigvee ^{k}_{i=1} \epsilon_{i,n}, \text{ with probability $p$}, \\
\left\{\bigvee ^{k}_{i=1}Y_{i,n-1}\right\} \bigvee \left\{\bigvee ^{k}_{i=1} \epsilon_{i,n}\right\}, \text{ with probability $(1-p)$.}
\end{cases}
\end{equation}

\noindent In terms of \textit{d.f}s and assuming stationarity this reads\\ \\
$F^k(x)=pF^k_\epsilon(x)+(1-p)F^k(x)F^k_\epsilon(x)$, or $F^k(x)=\frac{pF^k_\epsilon(x)}{1-(1-p)F^k_\epsilon(x)}$.\\That is, $F(x)=\left\{\frac{F^k_\epsilon(x)}{a-(a-1)F^k_\epsilon(x)}\right\}^{1/k}, a=\frac{1}{p}$. Hence we have,\\ 
\\
\textbf{Theorem 3.2} (Satheesh \textit{et al.} (2008)) A \textit{d.f} $F(x)$ can model the generalised stationary max-AR(1) scheme (3.1) for some $p\in(0,1)$ if it is Harris($a,k$)-maximum, $a=\frac{1}{p}$, and the distribution of the innovations is that of the components and conversely.\\
\\
\textbf{Theorem 3.3} (Satheesh \textit{et al.} (2008)) A \textit{d.f} $F(x)$ can model the generalised stationary max-AR(1) scheme (3.1) for every $p\in(0,1)$ (or as $p\downarrow 0$) if it is Harris($a,k$)-max-ID and conversely.\\
\\
\textbf{Remark 3.2} We saw that for $k>0$ integer and $a>0$, $H(x,\alpha) =  \left\{\frac{F^k(x)}{a-(a-1)F^k(x)}\right\}^{1/k}, x\in R$, is closed under Harris($b,k$)-maximum. Hence we can use the above to model the innovations in (3.1). The closure property implies also that if we repeat this operation with the same parameter $b>1$ they continue to be Harris($b,k$)-maximum and a passage to the limit shows that the limit is Harris-max-ID.\\

Now consider a variation of this max-AR(1) scheme described by \textit{i.i.d r.v}s $\{Y_{i,n}\}$ and innovations $\{\epsilon_{i,n}, i = 1,2,...,k\}$ for a fixed positive integer $k$ and some $c>0$.

\begin{equation}
\bigvee ^{k}_{i=1} Y_{i,n} = \begin{cases}
\bigvee ^{k}_{i=1}\frac{1}{c} Y_{i,n-1}, \text{ with probability $p$}, \\
\left\{\bigvee ^{k}_{i=1}\frac{1}{c}Y_{i,n-1}\right\} \bigvee  \left\{\bigvee ^{k}_{i=1} \epsilon_{i,n}\right\}, \text{ with probability $(1-p)$.}
\end{cases}
\end{equation}

In terms of \textit{d.f}s and assuming $\epsilon_{i,1}=Y_{i,0}$ and stationarity of $\{Y_{i,n}\}$, this reads\\

$F^k(x)=pF^k(cx)+(1-p)F^k(cx)F^k(x)$, or $F^k(x)=\frac{pF^k(cx)}{1-(1-p)F^k(cx)}.$\\
That is, $F(x)=\left\{\frac{F^k(cx)}{a-(a-1)F^k(cx)}\right\}^{1/k}, a=\frac{1}{p}$. Hence by induction we have\\ 
\\
\textbf{Theorem 3.4} For the stationary max-AR(1) model (3.2) to hold assuming $\epsilon_{i,1}=Y_{i,0}$, for some $p\in(0,1)$, the \textit{d.f} $F(x)$ (of $Y_{i,n}$) must be Harris($a,c,k$)-max-semi-stable and if we demand (3.2) is to be true for every $p\in(0,1)$, then $F(x)$ is Harris($a,k$)-max-stable and conversely.\\
\\
\textbf{Remark 3.3} By remark 2.2 we can construct Harris($a,c,k$)-max-semi-stable \textit{d.f}s.\\

\noindent{\textbf{Acknowledgement}}\\Authors thank the referee for asking to clarify a point that resulted in a better discussion in the paper. \vskip 0.5cm
\noindent{\textbf{References}}

\refmark Marshal, A W and Olkin, I (1997), On adding a parameter to a distribution with special reference to exponential and Weibull models, \textit{Biometrika}, \textbf{84}, 641-652.
\refmark Pancheva, E; Kolkovska, E T and Jordanova, P K (2006), Random time changed extremal processes, \textit{Probab. Theor. Appl.}, \textbf{51}, 752-772.
\refmark Pillai, R N (1990). Harmonic mixtures and geometric infinite divisibility, \textit{J. Ind. Statist. Assoc.}, \textbf{28}, 87-98.
\refmark Sandhya, E.; Sherly, S., Raju, N. (2008). Harris family of discrete distributions, in \textit{Some Recent Innovations in Statistics}, Published by: Department of Statistics, University of
Kerala, Trivandrum, India, 57–72.
\refmark Satheesh, S and Nair, N U (2004). On the stability of geometric extremes, \textit{J. Ind. Statist. Assoc.}, \textbf{42}, 99-109.
\refmark Satheesh, S; Sandhya, E and Rajasekharan, K E (2008), A generalization and extension of an autoregressive model, \textit{Statist. Probab. Lett.}, \textbf{78}, 1369-1374.

\end{document}